\documentclass[12pt,leqno]{amsart}

\usepackage{amsmath,stmaryrd}
\usepackage{amssymb}
\usepackage{amsthm} 
\usepackage{epsf}
\usepackage{graphicx}
\usepackage{esint} 
\usepackage{dsfont}
\usepackage[pagebackref]{hyperref}
\hypersetup{backref,  pdfpagemode=FullScreen,  colorlinks=true} 
\usepackage{color}
\usepackage{enumerate}

\numberwithin{equation}{section}

\newcommand{\R}{\mathbb{R}}
\renewcommand{\H}{\mathcal H}

\renewcommand{\d}{\partial}
\newcommand{\dist}{\,\mathrm{dist}}

\newcommand{\sm}{\setminus}

\newcommand{\diam}{\mathrm{diam}}

\newcommand{\wt}{\widetilde}

\newcommand{\ol}{\overline}

\newcommand{\cL}{{\mathcal  L}}

\newcommand{\cB}{{\mathcal  B}}

\newcommand{\1}{{\mathds 1}}

\newcommand{\dr}{\partial}
\DeclareMathOperator*{\osc}{osc}
\newcommand{\ds}{\displaystyle}

\theoremstyle{plain}
\newtheorem{theorem}[equation]{Theorem}

\theoremstyle{definition}

\theoremstyle{remark}

\newcommand{\re}{\mathbb{R}}

\newcommand{\RR}{{\mathbb{R}}}

\newcommand{\eps}{\varepsilon}

\newcommand{\A}{\mathcal{A}}
\newcommand{\C}{\mathcal{C}}

\newcommand{\po}{{\partial\Omega}}

\newcommand{\dint}{\int\!\!\!\int}

\DeclareMathOperator{\Tr}{Tr}
\DeclareMathOperator{\Ext}{Ext}
\def\div{\mathop{\operatorname{div}}}

\begin{document}

\title{A new elliptic measure on lower dimensional sets}
\author{
G. David, J. Feneuil, and S. Mayboroda}
\date{}
\begin{abstract}The recent years have seen a beautiful breakthrough culminating in a comprehensive understanding of certain scale-invariant properties of $n-1$ dimensional sets across analysis, geometric measure theory, and PDEs. The present paper surveys the first steps of a program recently launched by the authors and aimed at the new PDE approach to sets with lower dimensional boundaries. We define a suitable class of degenerate elliptic operators, explain our intuition, motivation, and goals, and present the first results regarding absolute continuity of the emerging elliptic measure with respect to the surface measure analogous to the classical theorems of C. Kenig and his collaborators in the case of co-dimension one.
\end{abstract}
\maketitle

\tableofcontents

\section{Introduction: history and motivation}

In the beginning of 80s 
 Carlos Kenig and his collaborators have launched a big program devoted to boundary value problems with data in $L^p$ and to a closely related issue of absolute continuity of harmonic measure with respect to the Hausdorff  measure on Lipschitz domains. It has been a bustling area of PDEs and harmonic analysis ever since, pushing the limits of allowable differential operators, geometry, function spaces. Somewhat later, roughly from the beginning of the 90s, 
  the developments at the interface of geometric measure theory  and harmonic analysis have been building up a deep understanding of scale-invariant properties of sets, and in particular, of the concept of uniform rectifiability. Finally, in the 21st century these somewhat independent pursuits combined and eventually culminated in a beautiful and powerful theory of equivalent scale-invariant characterizations of geometric, analytic, and PDE properties of sets. 

Unfortunately, most of these results are restricted to $n-1$ dimensional boundaries of $n$-dimensional domains and do not allow one to treat, for instance, a 1-dimensional curve in $\RR^3$. In 2012, 
 the authors of the present paper initiated a big project devoted to the properties of lower-dimensional sets in $\RR^n$. Here we survey some of our first results motivated by the aforementioned work of Kenig on the elliptic theory on Lipschitz domains. 

Let us start with the main brushstrokes of the results available in co-dimension 1, that is, for $n$-dimensional domains with $n-1$ dimensional boundaries. The pivotal work of Dahlberg \cite{Da77} addressing the absolute continuity of harmonic measure with respect to the Hausdorff measure on a Lipschitz domain brought a natural question of classifying all elliptic operators in the divergence form, $L=-{\rm div} A \nabla$, such that the corresponding elliptic measure has the same property. The main elements of the classical elliptic theory, Wiener's criterion, De Giorgi-Nash-Moser estimates, etc., suggest that all such operators behave similarly as long as the matrix of coefficients $A$ is bounded and elliptic. Absolute continuity, however, is a much more delicate property, and the counterexamples of Caffarelli, Fabes, Kenig \cite{CFK} and Modica, Mortola \cite{MM} showed that some regularity in $t$, the transversal direction to the boundary, is needed.  At that point the theory has split into two directions: the case of $t$-independent coefficients addressed for real symmetric matrices by Jerison and Kenig \cite{JK} and for real non-symmetric matrices by Hofmann, Kenig, Pipher, and the last of author of this paper in \cite{HKMP}, and the case of coefficients with controlled oscillations, satisfying the Carleson measure condition
\begin{equation}\label{CMcoeff}
\dint_{B(z, r)\cap \Omega}\sup\{\dist\{X, \po\}|\nabla A(X)|^2: \,X \in B(Z, \dist\{Z, \po\}/2)\} \,dZ \leq C r^{n-1},
\end{equation}
for all  $z\in \po,$ $r>0,$
pioneered by Kenig and Pipher in \cite{KePiDrift} (we will later provide more detailed references). Some of these results were further extended by Jerison and Kenig to a more general context on non-tangentially accessible domains \cite{JK2}. However, the true peak of the study of the domains with extremely rough geometry came much later, and now we  discuss some important recent achievements. In the present paper  we will primarily restrict our attention to quantitative results and to the realm of Ahlfors regular sets, and even in that setting we do not aim to provide a comprehensive review but rather some relevant highlights. We apologize in advance for many bibliographical omissions. 

A set $\Gamma \subset \R^n$ is a $d$-dimensional Ahlfors regular (AR) set if there exists a measure 
$\sigma$ supported on $\Gamma$ and a constant $C_0>0$ such that 
\begin{equation} \label{defADR}
C_0^{-1} r^d \leq \sigma(B(x,r)) \leq C_0 r^d
\ \text{ for $x\in \Gamma$ and $0 < r < \diam(\Gamma)$.} 
\end{equation}

When \eqref{defADR} holds for some measure $\sigma$, it is not hard to show that there
is a constant $C$, that depends on $n$ and $C_0$, such that 
$C^{-1}\H^d(A) \leq \sigma(A) \leq C \H^d(A)$ for $A \subset \Gamma$, where
$\H^d$ denotes the $d$-dimensional Hausdorff measure. Because of this, we may
replace $\sigma$ with the restriction $\H^d_{\vert \Gamma}$ of $\H^d$ to $\Gamma$, which also 
satisfies \eqref{defADR} for some $C_0' = C'_0(C_0,n)$. For simplicity we shall take 
$\sigma = \H^d_{\vert \Gamma}$, but all the results we shall talk about can be extended to any 
$\sigma$ satisfying \eqref{defADR}.

As mentioned above, for Ahlfors regular sets of co-dimension 1 (i.e. of dimension $n-1$), there are precise results that prove the equivalence of scale-invariant geometric properties of $\Gamma$, 
a certain form of the harmonic analysis on $\Gamma$, and the good behavior of partial differential equations on $\R^n \sm \Gamma$. 
To be specific, the equivalent geometric properties that 
\begin{enumerate}[(G1)]
\item $\Gamma$ is uniformly rectifiable (see \cite{DS1,DS2}, but note that similar properties were considered before, starting with the Lipschitz, chord-arc, or big pieces conditions of \cite{David88,David88b,David91}), 
\item P. Jones $\beta$-numbers of \cite{Jones} 
satisfy a suitable form of a Carleson measure condition \cite{Jo2,Ok,DS1,DS2}, 
\item Wasserstein distance functions $\alpha$ of X. Tolsa satisfy a suitable Carleson measure condition \cite{Tolsa09}, 
\end{enumerate}
hold if and only if any of the following analytic properties is valid 
\begin{enumerate}[(H1)]
\item the $L^2$-boundedness on $\Gamma$ of all the standard odd Calder\'on-Zygmund singular integral operators,
\item more specifically the $L^2$-boundedness on $\Gamma$ of the Riesz transform alone,
\item usual square function estimates for the Cauchy/Newtonian kernel. 
\end{enumerate}
The fact that uniform rectifiability implies the conditions above has been known for some time
\cite{CMM,David88,David88b,David91,Se,DS1}, as well as the equivalence of (G1)--(G3) and (H1), (H3)
\cite{DS1,Tolsa09}. More recently, the implication (H2) $\Longrightarrow$ (H1) was proved (see \cite{MMV} in the complex plane and \cite{NToV} in higher dimensions) which opened a door to the PDE properties that we discuss below.

For the present paper we are more interested in the connection of all these with
the regularity properties of elliptic PDE's on $\Omega = \R^n \sm \Gamma$. The following are equivalent to (G1)--(G3) and (H1)--(H3), in particular, 
to the uniform rectifiability of $\Gamma$: 
\begin{enumerate}[(P1)]
\item all bounded harmonic functions in $\Omega$ satisfy Carleson measure estimates \cite{GMT}, \cite{HMM1},
\item all bounded harmonic functions in $\Omega$  are $\eps$-approximable \cite{GMT}, \cite{HMM1}.
\end{enumerate}

We mention above the papers which established the final equivalence results, however, we would like to underline the importance of the work by Kenig, Kirchheim, Koch, Pipher, and Toro \cite{KKiPT}, \cite{KKPT}, which established equivalence of these conditions to the absolute continuity of harmonic measure with respect to the Lebesgue measure on a Lipschitz graph, both motivating the attention to properties (P1)--(P2) on more general domains and providing the essential tools for their analysis. We mention, parenthetically, that we should rather be talking of an appropriate quantifiable analogue of absolute continuity, the $A^\infty$ condition, but let us not dive into such details for the purposes of the introduction.

However, in the general context of Ahlfors regular sets or even that of uniformly rectifiable domains, the absolute continuity of harmonic measure with respect to the Hausdorff measure on the boundary is not equivalent to (P1), (P2) or to any other condition above, and, in fact, may fail in a complementary component $\Omega$ of a uniformly rectifiable set \cite{BJ}. In short, one needs some form of quantifiable connectedness of $\Omega$ to ensure reasonable access to the boundary. A search for the full, comprehensive understanding of this elusive topological condition has motivated many beautiful results, as well as the study of (P1)--(P2) in the first place, but we will only mention the final characterizations that emerged just this year: given an $(n-1)$ Ahlfors regular set $\Gamma$ and $\Omega = \R^n \sm \Gamma$, the harmonic measure is  $A^\infty$ on $\Gamma$ if and only if the boundary $\Gamma$ is uniformly rectifiable and $\Omega$ is semi-uniform  \cite{Azz}; the harmonic measure is locally weak-$A^\infty$ on $\Gamma$ if and only if the boundary $\Gamma$ is uniformly rectifiable and $\Omega$ satisfies the weak local John condition  \cite{HM18, AzMouTo}. Finally, inspired by the aforementioned earlier achievements of Kenig and Pipher \cite{KePiDrift}, many of these results have been further extended to operators $L=-\div  A \nabla$ with the coefficients satisfying the Carleson measure condition \eqref{CMcoeff}  (notice that $t$-independence hardly makes sense in such a general geometric setting).  We shall shortly return to this in Section \ref{S15}, and, in particular, give some of the relevant definitions. For now let us move to the subject of this paper. 

\vskip 0.08 in

Virtually none of the aforementioned characterizations has been available for sets of 
co-dimension higher than one. The main  reason for this is that the standard elliptic
operators, such as the Laplacian, do not interact well with higher-dimensional boundaries $\Gamma$, and if fact, even seemingly unrelated results regarding (H3) rely heavily on the harmonicity of the Riesz transform.
The Dirichlet problem, for instance, is hard to solve in such a context, since reasonably sized 
harmonic functions on $\Omega = \R^n \sm \Gamma$ have a tendency to extend across $\Gamma$ because 
it is too small, or, said differently, brownian paths almost never meet $\Gamma$. 

Lewis, Nystr\"om, and Vogel \cite{LN,LNV}
have found a nice way to tackle this problem by replacing $\Delta$ with a $p$-Laplacian,
but fine relations between the geometry and the operator seem to be very hard to come
with in this context. In particular, it is not clear at this point whether the resulting $p$-harmonic measure is absolutely continuous with respect to the Hausdorff measure even on a Lipschitz graph of dimension $d<n-1$. We propose a different approach.

A few years ago the authors of the present paper have launched a theory of degenerate PDEs 
on sets with lower dimensional boundaries, which 
gives hope to extend some of the aforementioned results to the setting of higher co-dimension, including, in some instances, fractal boundaries and fractional dimension scenarios. 
On top of the considerably generalized elliptic theory and expected (but now much more involved) results on lower-dimensional Lipschitz graphs, our approach
brings out surprising new theorems in the ``classical" setting of co-dimension 1 and even
introduces some new mysterious identities. 
It merges, at some level, with the theory of non-local operators and the fractional Laplacian 
which have recently attracted a lot of attention
albeit in a completely different context, and intertwines in various ways with the results of Kenig and Pipher \cite{KePiDrift}, Kenig, Kirchheim, Pipher, Toro \cite{KKPT},  as well as Fabes, Kenig, Serapioni \cite{FKS}. In some sense, consideration of a full class of operators satisfying a suitable modification of the Carleson measure condition \eqref{CMcoeff} now becomes a necessity even for the simplest foundational results. At present, we feel that we have only scratched the surface of a vast  field of investigation. 
 
This  survey will discuss the first results and  is 
based on the work on the three authors, together with two of their collaborators: 
Zihui Zhao and Max Engelstein. 

In Section \ref{S15} we consider boundaries of co-dimension $1$ for the last time and rapidly review some of the elliptic theory that we want to generalize. 

Section \ref{S2} presents the basic set-up of the elliptic theory in our context.
We discuss the homogeneity and other reasonings behind our choice of operators, the construction of 
the elliptic measures in the general context of Ahlfors regular sets (possibly of fractional dimension), and the basic properties of these elliptic measures. 
This is based on \cite{DFMprelim}, except for some small additional properties that 
come from 
\cite{DFMprelim2}. We believe that this is the right 
setting for studying the relations between PDEs and geometry presented in the introduction, 
that hold for sets of co-dimension 1.

Section \ref{S3} is devoted to the solvability of the Dirichlet problem in $\R^n \setminus \Gamma$ with data in $L^p$,  and we will rapidly concentrate on the special case when $\Gamma \subset \R^n$ is a $d$-plane,
and the matrix $A$ associated to $L=-\div A \nabla$ satisfies some appropriate Carleson
measure condition. The results presented in this section (Theorems \ref{3.12} and \ref{3.16})
generalize results (of co-dimension $1$) for elliptic operators in the upper half space that seem
to be new too. They are the main aim of \cite{FMZ}, which heavily uses \cite{DFMAinfty} (discussed again below) and a perturbation theory similar to the one developed in \cite{DPP} in co-dimension 1.

In Section \ref{S4}, we return to more general Ahlfors regular boundaries $\Gamma$ of dimension $d < n-1$,
and explain that when $\Gamma$ is the graph of Lipschitz function $\varphi : R^d \to \R^{n-d}$ with
small enough Lipschitz norm, and $L$ is a rather natural degenerate elliptic operator defined
in terms of $\Gamma$, the elliptic measure of Section \ref{S2} is (mutually) absolutely continuous 
with respect to $\sigma = \H^d_{\vert \Gamma}$, with an $A_\infty$ weight. 
This can be seen as the analogue of the result of Dahlberg \cite{Da77} for sets with lower dimensional boundaries, but notice that even the proper analogue of the Laplacian is not automatically clear in our context. The proof \cite{DFMAinfty} relies on a new change of variables that sends
$\Omega = \R^n \sm \Gamma$ back to $\Omega_0 = \R^n \sm \R^d$ virtually conformally, or in other words, carefully preserving an isometry in the transversal direction to the boundary. We also talk quickly about the Dirichlet boundary problem with data in BMO (\cite{MZ}).

At last, Section \ref{S6} briefly discusses some first analogues of the characterizations (G1)--(H4) that we are considering.

\section{Elliptic operators and boundaries of co-dimension $1$}
\label{S15}

Let us recall here the standard setting for elliptic operators on a domain $\Omega \subset \R^n$
that we want to generalize. For the simplicity of the exposition, we shall make the following geometric
assumptions on $\Omega$. First, that
\begin{equation} \label{EAR}
\text{$\Gamma = \d\Omega$ is Ahlfors regular of dimension $n-1$,}
\end{equation}
which means that \eqref{defADR} holds, say, for $\sigma = \H^{n-1}_{\vert \Gamma}$.
This is a strong quantitative way to require that $\Gamma$ is $(n-1)$-dimensional. 
This may be too strong, but then statements are simpler and certainly if we want
to compare the harmonic measure to $\sigma$, it makes sense to assume that $\sigma$
is not too bad. 

We also assume that both $\Omega$ and $\R^n \sm \ol\Omega$ satisfy the corkscrew
condition. We say that $\Omega$ satisfies the corkscrew point condition 
(quantitative openness) when there exists a constant $C_1 > 0$ such that
for $x\in \d\Omega$ and $0 < r < \diam(\d\Omega)$, 
\begin{equation} \label{2.7}
\text{one can find a point $A_{x,r} \in \Omega \cap B(x,r)$ such that 
$B(A_{x,r},C_1^{-1} r) \subset \Omega$.}
\end{equation}
Finally, we require that $\Omega$  satisfies the Harnack chain condition (quantitative connectedness).
This means (also see \cite[Page 13, property (iii)]{KenigB}) 
that there is a constant $C_2 \geq 1$ and, for each $\Lambda \geq 1$, 
an integer $N \geq 1$ such that, whenever $X,Y\in \Omega$ and $r \in (0, \diam(\d\Omega))$ 
are such that
\begin{equation} \label{cs1} 
\min\{\dist(X,\d\Omega),\dist(Y,\d\Omega)\} \geq r \ \text{ and } \  |X-Y| \leq \Lambda r,
\end{equation}
we can find a chain of $N+1$ points $Z_0=X, Z_1,\dots, Z_{N} = Y$ in $\Omega$ such that
\begin{equation} \label{cs2}
C_2^{-1} r \leq \dist(Z_i,\d\Omega) \leq C_2 \Lambda r
\text{ and } 
|Z_{i+1}-Z_i| \leq \dist(Z_i,\d\Omega)/2 
\text{ for } 1 \leq i \leq N.
\end{equation}
These assumptions are too strong compared to the modern state of the art as described in the introduction, but  they nicely parallel some consideration we will have for the sets with lower dimensional boundaries, so let us suppose for now that all of the geometric conditions above are satisfied in $\Omega$.

In addition to the harmonic measure associated to the Laplacian, we want to study
the elliptic measures associated to operators $L=-\div  A  \nabla$, and we shall
restrict our attention to real-valued, uniformly elliptic matrices  $A$ with bounded measurable coefficients. To be specific, we assume
that $A$ is measurable on $\Omega$, takes values in the set of $n \times n$ matrices,
and that there is a constant $C_3 \geq 1$ such that
\begin{equation} \label{2.19}
\big|A(X) \xi \cdot \zeta \big| \leq C_3 |\xi||\zeta| 
\quad \text{ for $X\in \Omega$ and $\xi,\zeta \in \R^n$}
\end{equation}
and
\begin{equation} \label{2.20}
A(X) \xi \cdot \xi \geq C_3^{-1} |\xi|^2 \quad \text{ for $X\in \Omega$ and $\xi \in \R^n$.}
\end{equation}

Let us first say a few words about a standard way to define the elliptic measure associated
to $L$ in this case. Further assume, just for the definition below and to make things even simpler, 
that $\Omega$ is also bounded.
Under these conditions, it is well known (see \cite{Stampacchia65}) that for any 
$g\in C^\infty(\partial \Omega)$, we can find a unique 
$u \in W^{1,2}(\Omega) \cap C(\overline{\Omega})$ such that 
\begin{equation} \label{2.1}
\left\{ \begin{array}{ll} 
Lu = 0 & \text{ in } \Omega \\
 u = g & \text{ on } \dr \Omega.
\end{array} \right.
\end{equation}
By a combination of the maximum principle and the Riesz representation theorem, one infers that 
for $X \in \Omega$, there is a 
Borel regular probability measure $\omega^X_L$ 
on $\d\Omega$ such that the solution of \eqref{2.1} satisfies
\begin{equation} \label{2.2}
u(X) = \int_{\dr \Omega} g(y) \, d\omega^X_L(y).
\end{equation}
The measure $\omega^X_L$ is called the elliptic measure with pole at $X$; 
the dependence of $\omega^X_L$ on $X$ is secondary, since the Harnack inequality implies
that $\omega^X_L$ are equivalent for two different choices of $X$. 
Note also that $u_E(X):= \omega_L^X(E)$ can be seen formally as the solution of \eqref{2.1} for $g:= \1_E$. The construction and the main
properties of the elliptic measures can be found for instance in C. Kenig's book \cite[Sections 1.1 to 1.3]{KenigB}. 

Under the geometric conditions above the elliptic measure is well-defined for any bounded elliptic $A$, the solutions satisfy the boundary and interior De Giorgi-Nash-Moser estimates, that is, local quantitative boundedness and H\"older continuity properties, the comparison principle holds, as well as other classical properties.   

The natural notion of quantitative mutual absolute continuity in the context of elliptic measures on AR sets is the following.
We shall say that $\omega_L$ and $\sigma$ are $A_\infty$-absolutely continuous
(in short, that $\omega_L \in A_\infty(\sigma)$) when for every 
$\epsilon>0$, there exists $\delta>0$ such that for any ball $B(x,r)$ centered
on $\Gamma$ and with radius $r\in (0,\diam(\Gamma))$, and any measurable set $E \subset \Gamma \cap B(x,r)$,
one has
\begin{equation} \label{3.23}
\omega_X^{A_{x,r}}(E) < \delta \Longrightarrow \frac{|E|}{|B|} < \epsilon,
\end{equation}
where $A_{x,r}$ is the corkscrew point of \eqref{2.7}. It can be checked that this property
does not depend on the choice of $A_{x,r}$.

In the particular case where $L = - \Delta$, the elliptic measure is called the harmonic measure.
Under the geometric assumptions above, or even less, assuming AR, {\it interior} corkscrews, and Harnack chains,  
 $\omega_\Delta \in A_\infty(\sigma)$ if and only if $\Gamma$ is
uniformly rectifiable if and only if $\Omega$ is chord-arc (i.e., in fact, it satisfies the exterior corkscrew condition as well)  \cite{AHMNT}. This is an anterior result to those we have listed in the introduction, but again, a good guidance for what we are planning to do next. 
But we may also be interested in other elliptic operators $L$, and in particular in the conditions 
that we need to put on the matrix $A$ if we want the condition $\omega_L \in A_\infty(\sigma)$
to follow from, or even be equivalent to, uniform rectifiability.  As mentioned in the introduction, in view of the counterexamples \cite{CFK}, \cite{MM} it makes sense to impose some regularity in the transversal direction to the boundary, and while the concept of $t$-independence does not make much sense in the absence of a special transversal direction, the Carleson measure conditions \eqref{CMcoeff} generalize to this setting and should be the right background hypotheses. And indeed, combining the results of Kenig, Pipher \cite{KePiDrift} with \cite{DJ} and \cite{AHMNT}, we conclude that the elliptic measure  $\omega_L \in A_\infty(\sigma)$ if the operator satisfies \eqref{CMcoeff}. The converse, currently only under the assumption of the smallness of Carleson measure in \eqref{CMcoeff}, has been established in
\cite{HMMTZ}. This is a somewhat rough outline, omitting certain technical details, e.g., on has to assume that $A\in {\rm Lip}_{loc}(\Omega)$ and $|\nabla A| \dist(\cdot, \po)\in L^\infty(\Omega)$, but since a thorough survey of the codimension 1 theory is not our goal, we remain a little imprecise at this point. Starting with the next chapter, the statements will become complete and detailed.

\section{Elliptic measure on sets with co-dimension higher than 1}
\label{S2}

\subsection{A weight $w$ with the right homogeneity}
\label{S2.1}

If we want to work with the Laplacian or standard elliptic operators, we will not be able to
discuss harmonic measure for AR sets of high co-dimensions. Let us rapidly say why, and then
how we intend to fix this problem.

Consider for example $\Gamma_0 = \{(x,0)\in \R^n, \, x\in \R^d\}$, with $d < n-1$.
First try to use the Laplacian on $\Omega = \Omega_0 = \R^n \sm \Gamma_0$.
This is not a bounded domain any more, but we can use homogeneous spaces
to deal with this. As customary, we understand the equation $Lu= - \Delta u = 0$ in \eqref{2.1}
in the weak sense, which means that  
\begin{equation} \label{2.3}
\int_\Omega \nabla u \cdot \nabla \varphi \, dX = 0 \qquad \text{ for }\varphi \in C^\infty_0(\Omega),
\end{equation}
and the reasonable space for doing this is 
$W_\Omega := \big\{u\in L^1_{loc}(\Omega),\, \nabla u \in L^2(\Omega)\big\}$.
But our boundary $\Gamma_0$ is now so small that all the functions of $W_\Omega$
extend through $\Gamma_0$, so $W_\Omega$ is equal to the {\em a priori} smaller space 
$W_{\R^n}=\big\{u\in L^1_{loc}(\R^n),\, \nabla u \in L^2(\R^n)\big\}$, and that
any function in $C^\infty_0(\R^n)$ can be approached, in term of the semi-norm 
$\|\nabla \cdot\|_2$, by functions in $C^\infty_0(\Omega)$. Then \eqref{2.3} implies that
\begin{equation} \label{2.4}
\int_\Omega \nabla u \cdot \nabla \varphi \, dX = 0 \qquad \text{ for }\varphi \in C^\infty_0(\R^n)
\end{equation}
and therefore the only solutions $u \in W_\Omega$ that satisfy \eqref{2.4} are harmonic functions
in $W_{\R^n}$, hence constant functions. In short, 
the boundary is too thin, \eqref{2.1} has no non-constant solutions, and therefore we cannot 
build the harmonic measure (or similarly the elliptic measure) as we did before. 

Another way to put this is that functions of $W_\Omega$ do not have well defined traces
on $\Gamma_0$, which makes it harder to define a Dirichlet problem; working in a smaller space than
$W_\Omega$ could lead to functions with a trace, but then we could only work with Dirichlet
data that are the restriction of harmonic functions near $\Gamma_0$, which is not what we intended.
Let us finally mention that the pleasant definition of the harmonic measure of $A \subset \Gamma_0$
as the probability that a Brownian path starting from $X$ hits $A$ the first time it leaves $\Omega$
does not work either, because it is well known that the Brownian path will almost surely never
leave $\Omega$.

\medskip

The strategy of the authors is to replace the classical elliptic operators by degenerate elliptic operators, also of the form $L=-\div  A  \nabla$, but where $A$ satisfies ellipticity and boundedness 
condition with a different homogeneity. That is, we shall define a weight $w(X)$ in a moment
(see \eqref{defw}), depending on the geometry of $\Omega$, and we require that for some 
$C_3 \geq 1$
\begin{equation} \label{2.5}
A(X) \xi \cdot \zeta \leq C_3 w(X) |\xi||\zeta| \quad \text{ for $X\in \Omega$ and $\xi,\zeta \in \R^n$}
\end{equation}
and
\begin{equation} \label{2.6}
A(X) \xi \cdot \xi \geq C_3^{-1} w(X) |\xi|^2 \quad \text{ for $X\in \Omega$ and $\xi \in \R^n$.}
\end{equation}
When $\Omega \subset \R^n$ has a $(n-1)$-dimensional Ahlfors regular boundary, 
we choose $w(X) \equiv 1$ and we recover the classical theory. 
But we shall work with $\Omega = \R^n \sm \Gamma$ for some AR set $\Gamma$ of dimension $d < n-1$,
and we will choose a weight $w(X)
$ that goes to $+\infty$ as $X$ gets close to the boundary. 
The effect of this is to increase the ``visibility'' of the boundary $\Gamma$, both in terms
of the associated function spaces (as we shall see soon) and, it seems to the authors, 
in terms of the stochastic process that replaces the Brownian motion, which will tend to 
be attracted in the direction of $\Gamma$. 

Of course we expect the precise speed at which $w(X)$ tends to $+\infty$ to be important,
and we shall now explain how we could guess a good choice for $w$ by looking at our special example 
where $\Gamma = \Gamma_0  = \{(x,0)\in \R^n, \, x\in \R^d\}$ and $\Omega = \R^n \sm \Gamma_0$
for some integer $d < n-1$. Let us even focus on radial solutions, i.e., solutions of \eqref{2.1}
that can be written $u(x,t) = u_0(x,|t|)$, where $u_0$ is defined on $\R^{d+1}_+$ 
and satisfies $u_0(x,0) = u(x,0) = g(x)$. 

Now $u_0$ has to satisfy some equation, if possible elliptic, and the simplest is to 
look for an example for which \eqref{2.1} for radial functions becomes 
$-\Delta u_0 = 0$ in $\R^{d+1}_+$. A computation with cylindrical coordinates
shows that we can take $L_0 := - \div |t|^{d+1-n} \nabla$. That is, we have a preferred
``substitute'' for the Laplacian, which is the operator $L_0$, and the computations
show that radial solutions of \eqref{2.1} come from solutions of $-\Delta u_0 = 0$ in $\R^{d+1}_+$
such that $u_0(x,0) = g(x)$. This suggests that we take 
\begin{equation} \label{defw}
w(X) = \dist(X,\d\Omega)^{d+1-n},
\end{equation}
at least in this case.

The reader will probably agree that the most beautiful operator to take in our special
domain is $L_0$; for more general domains $\Omega = \R^n \sm \Gamma$, where $\Gamma$ is
AR of dimension $d < n-1$, we often find it hard to select a preferred operator like $L_0$ (for instant, to start a study of absolute continuity, much as one typically starts with the Laplacian in co-dimension 1), 
even though we shall try to present some that we like best in the next sections.
However, first, let us try to establish the analogues of the maximum principle, local regularity of solutions, definition of elliptic measure, and other fundamental properties of elliptic theory, for the more general class of elliptic operators, taking $w(X) = \dist(X,\d\Omega)^{d+1-n}$ as in \eqref{defw}, and requiring that
$A$ satisfy \eqref{2.5} and \eqref{2.6} with that weight.

\subsection{Basic elliptic theory I: our assumptions, the Geometry of $\Omega$, the two relevant Hilbert spaces, traces, extensions, definitions of solutions}

Once we have the weight $w$ and the corresponding class of elliptic operators,
we can try to follow some classic route to define elliptic measure and study its basic properties.
The theory presented in the next two paragraphs, which is a summary of \cite{DFMprelim,DFMprelim2}, 
has some strong connection with the work of Fabes, Jerison, Kenig, Serapioni in the beginning 
of the 1980's (see \cite{FKS,FJK,FJK2}) who have laid some foundations of the study of degenerate elliptic operators. 
However, even though there are some overlaps between the two works, the spirit is  
different, in that Fabes and others studied to which extend the elliptic theory can be recovered if one is given a degenerate elliptic operator (in particular they gave a Wiener test), while we choose a weight $w$ (and thus the degeneracy of the elliptic operator) adapted to the boundary so that, among other properties, the domain always satisfies the Wiener test.

We fix for the rest of the section a dimension $d<n-1$, which is not necessarily an integer. 
Let $\Gamma$ be a $d$-dimensional Ahlfors regular set, as in \eqref{EAR}, and consider  
the domain $\Omega = \R^n \setminus \Gamma$. We assume $\Gamma$ to be unbounded, so
as to fit with \cite{DFMprelim,DFMprelim2}.
Set 
\begin{equation} \label{delta}
\delta(X) = \dist(X,\Gamma) \ \text{ for } X \in \Omega,
\end{equation}
and then define the weight $w$  
by $w(X) = \delta(X)^{d+1-n}$ as in \eqref{defw}.

In this subsection we take care of some simple geometry,
define the two associated Hilbert spaces $W$ (our energy space) and $H$ (the space of traces),
and operators (of trace and extension) between these spaces. 

The thinness of 
$\Gamma$ has advantages, in that there is plenty of space inside the domain. 
In particular, we shall get for free that $\Omega$ satisfies the corkscrew point condition 
\eqref{2.7}, and the following Harnack chain condition: 
for any choice of $X,Y\in \Omega$ and $r>0$ such that $\min\{\delta(x),\delta(Y)\} \geq r$ 
and $|X-Y| \leq 7C_1 r$, where $C_1$ comes from \eqref{2.7}, 
we can find a chain of $C_2 +1$ points $Z_0, \dots, Z_{C_2}$ in $\Omega$, such that
\begin{equation} \label{2.8}
Z_0 = X, Z_{C_2} = Y, \ \text{ and $|Z_{i+1}-Z_i| \leq \delta(Z_i)/2$ for $0 \leq i < C_2$.}
\end{equation}
Here $C_1$ and $C_2$ depend only on $C_0$ and $n-d-1$.

The proof of \eqref{2.7} is done by using the AR property \eqref{EAR} to estimate the maximum number
of disjoint balls of radius $\varepsilon r$ that meet $\Gamma \cap B(x,r)$; then we choose 
$\epsilon := \epsilon (C_0,n-d-1)$ so small that at least one such ball contained in $B(x,r)$
does not meet $\Gamma \cap B(x,r)$.
As for the proof of \eqref{2.8}, we show similarly that we can find a tube $T$ of radius $\epsilon r$, 
$\epsilon:=\epsilon(C_0,n-d-1,C_1)$ small enough, that links $B(X,r/2)$ to $B(Y,r/2)$ without meeting
$\Gamma$; the points $Z_i$ are then taken appropriately in $B(X,r) \cup B(Y,r) \cup T$.

Let us comment more on the Harnack chain condition. The points of \eqref{2.8} yield a 
chain of balls $B_j = B(Z_i,\delta(Z_i)/2)$, $0\leq i \leq C_2$, called a Harnack chain,
and we easily observe that the chain stays at distance at least $2^{-C_2}r$ from $\Gamma$,
and does not go too far either. At first sight, what we get with \eqref{2.8} seems weaker that
what we announced in \eqref{cs1} and \eqref{cs2}, because we do not allow $X, Y$ such that
$|X-Y| \leq \Lambda r$ for $\Lambda$ larger than $7C_1$. But this is not a problem because
we can iterate the construction. In fact, if $X$ and $Y$ satisfy \eqref{cs1},
one can find at most $C'_2 \ln(1+\Lambda) +1$ points $W_0, \dots, W_{N}$ in $\Omega$,
with $W_0 = X$, $W_N = Y$, such that 
\begin{equation} \label{2.9}
\text{$|W_{i+1}-W_i| \leq \delta(W_i)/2$ and 
$2^{-C_2}r \leq \delta(W_i) \leq C_1 2^{C_2+4}\Lambda r$ for $0 \leq i < N$,}
\end{equation}
which is better than \eqref{cs2}.

Let us sketch the proof of this fact.  Set $x,y\in \Gamma$ be such that $\delta(X) = |X-x|$ and $\delta(Y) = |Y-y|$.  Thanks to \eqref{2.7}, we construct $X_i,Y_j \in \Omega$ as a corkscrew point associated to  $(x,2^i C_1 \delta(X))$ and $(y,2^j C_1 \delta(Y))$ respectively.  We choose $i_0,j_0$ such that $2^{i_0}\delta(X) \approx 2^{j_0}\delta(X) \approx \Lambda$.  Two successive points of the sequence $X,X_0,\dots, X_{i_0},Y_{j_0},\dots,Y_0,Y$ -- which has a length smaller than $C\ln(1+\Lambda)+1$ -- are actually close enough to each other to allow the use of the chain property with \eqref{2.8} and construct a small Harnack chain linking them.  The concatenation of these small chains gives a long chain that satisfies \eqref{2.9}.

In summary, if $\Gamma$ is AR of dimension $d < n-1$, then
$\Omega = \R^n \sm \Gamma$ is 
a uniform domain, which means that it satisfies
the corkscrew condition \eqref{2.7} and the Harnack chain condition. 

As we discussed above, in co-dimension $1$ it is important for the complement of the domain that we consider to be fat enough, e.g., that the {\it exterior} corkscrew condition is satisfied (in fact, such an exterior condition is even necessary for some PDE properties). With such an additional requirement, $\Omega$ is commonly called Non Tangentially Accessible (NTA), and NTA domains offer a rather friendly setting for the elliptic PDEs. However, our $\Omega$ does not satisfy any exterior fatness condition, since  $\R^n \sm \Omega = \Gamma$ is much too small, so we need to compensate.

Our salvation will come from the fact that we are able to avoid, or rather reformulate, a condition of fatness of the complement of the domain by proving a boundary Poincar\'e inequality, which in some sense establishes that a suitable capacity of the boundary must be massive. 
More precisely, we shall prove that for any ball $B=B(x,r)$ centered on $\Gamma$ 
and any smooth function $u$ that takes the value $0$ on $\Gamma \cap B$,
one has a Poincar\'e inequality of the form
\begin{equation} \label{2.10}
\int_{B} |u|^2 w(X) dX 
\leq C r^2 \int_B |\nabla u|^2 w(X) dX.
\end{equation}
This Poincar\'e estimate would fail for $w\equiv 1$ for the same reason as \eqref{2.1} 
cannot be solved with $L = -\Delta$, and our present choice of $w = \delta(X)^{d+1-n}$ 
also comes from the fact that we were actually able to prove \eqref{2.10} with such weights $w$. 

As suggested by \eqref{2.10}, the most relevant measure on $\Omega$ is
defined by $dm(X) = w(X)dX$, or equivalently by
\begin{equation} \label{mdX}
m(E) := \int_E w(X)dX \ \text{ for } E \subset \Omega.
\end{equation}
It is not very hard to see that $m$ is a doubling measure on $\R^n$, and we even have 
\begin{equation} \label{2.12}
C^{-1} \left( \frac{r}{s} \right)^{d+1} \leq \frac{m(B(X,r))}{m(B(X,s))} \leq C \left( \frac{r}{s} \right)^{n} \qquad \text{ for $X\in \R^n$ and $0<s<r$,}
\end{equation}
and $m(B(x,r)) \approx r^{d+1}$ for $x\in \Gamma$ and $r>0$. 

We now come to the definition of the energy space $W$. 
Since in the higher co-dimension case we systematically work with an unbounded domain $\Omega$, we shall use homogeneous spaces. To this end, define the weighted Sobolev space $W$ by
\begin{equation} \label{2.13}
W:= \big\{ u \in L^1_{loc}(\Omega), \, \nabla u \in L^2(\Omega,dm)\big\},
\end{equation}
which happens 
to be equal to the apparently smaller space 
$\big\{ u \in L^2_{loc}(\R^n), \, \nabla u \in L^2(\R^n,dm)\big\}$. 
Set $\| u \|_W = \|\nabla u\|_{L^2(dm)}$; 
the quotient of the semi-normed space $(W,\|.\|_W)$ by constant functions is a Hilbert space.  
In addition, due to the fact that $W \subset L^1_{loc}(\R^n)$, 
we can use mollifiers to prove that 
\begin{equation} \label{2.14}
\begin{split}
& \text{any function in $W$ can be approached in the semi-norm $\|.\|_W$} \\
& \hskip2cm \text{by functions in $C^\infty(\R^n) \cap W$.}
\end{split}
\end{equation}

Our second Hilbert space is the set of traces $H$, which is the set of measurable
functions $g$ defined on $\Gamma$ and such that  
\begin{equation} \label{2.15}
\|g\|_H^2 
:= \int_\Gamma\int_\Gamma \frac{|g(x)-g(y)|^2}{|x-y|^{d+1}} d\sigma(x) \, d\sigma(y) < +\infty ,
\end{equation}
where we recall that $\sigma=\mathcal H^d_{\vert\Gamma}$, or any measure satisfying \eqref{defADR}.
The expert reader will have recognized a definition of the Sobolev space $H = H^{1/2}(\Gamma)$
associated to the measured metric space $(\Gamma, \dist, \sigma)$.

From here on, the functional analysis will go pretty smoothly.
We can construct two bounded linear operators, the trace operator $\Tr: \, W\to H$ and the extension operator $\Ext: \, H\to W$ such that $\Tr \circ \Ext = {\rm Id}_H$. The trace of $u\in W$ is such that for $\sigma$-almost every $x\in \Gamma$, 
\begin{equation} \label{2.16}
\Tr u(x) = \lim_{r\to 0} \fint_{B(x,r)} u(Y) dY,
\end{equation}
and even, analogously to the Lebesgue density property,
\begin{equation} \label{2.17}
\lim_{r\to 0} \fint_{B(x,r)} |u(Y) - \Tr u(x)| dY = 0.
\end{equation}
The completion of $C^\infty_0(\Omega)$ with the norm $\|.\|_W$ is the space 
$W_0 := \big\{f \in W, \, \Tr f = 0 \big\}$. Using results from \cite{HK2,HK00}, 
the Poincar\'e inequality \eqref{2.10} can be improved into a Sobolev-Poincar\'e inequality that 
holds for any function $f\in W$, which means that 
\begin{equation} \label{2.18}
\left(\fint_{B(x,r)} |u|^p dm\right)^{1/p} \leq C_p r \left(\fint_{B(x,r)} |\nabla u|^2 dm\right)^\frac12
\end{equation}
for $x\in \Gamma$, $r>0$, $u\in W$ such that $\Tr u = 0$ on $\Gamma \cap B(x,r)$, and $p\in [1,\frac{2n}{n-2}]$ (if $n\geq 3$) or $p\in [1,+\infty)$ (if $n=2$). The property \eqref{2.14} allows us to get a chain rule: let $f\in C^1(\R)$ be such that $f'$ is uniformly bounded on $\R$ and let $u \in W$, then $f\circ u \in W$ and $\nabla (f\circ u) = f'(u) \nabla u$ a.e. in $\R^n$. This 
is a key property that allows us to consider the quantities $|u|$ or $\max\{k,u\}$, needed in the proof of the De Giorgi-Nash-Moser estimates we shall see in the next paragraph.  

\subsection{Basic elliptic theory II: solutions, elliptic measure, and their properties}

We keep the same assumptions and notation as before for $\Gamma$ 
(an AR set of dimension $d < n-1$), the weight $w$ defined by \eqref{defw}, and the two 
associated Hilbert spaces $W$ and $H$. The geometry is under control, and we now consider 
any degenerate elliptic operator $L = -\div A \nabla$, with measurable real matrix-valued 
coefficients $A$ that satisfy \eqref{2.5} and \eqref{2.6}. 
Often we prefer to use the rescaled elliptic matrix $\A = w^{-1}A$, 
which therefore satisfies the classical conditions \eqref{2.19} and \eqref{2.20}.

A function $u\in W$ is called a weak solution to $Lu = 0$ in $\Omega$ if
\begin{equation} \label{2.21}
\int_\Omega A \nabla u \cdot \nabla \varphi \, dX = \int_\Omega \A \nabla u \cdot \nabla \varphi \, dm  = 0 \qquad \text{ for } \varphi \in C^\infty_0(\Omega),
\end{equation}
where the first equality comes directly from the definition of $\A$. 
By \eqref{2.19} and \eqref{2.20} for $\A$, it is clear that $W$ is the right space for the use of 
the Lax-Milgram theorem. Using the latter and the extension operator $\Ext$, we see that
for each $g\in H$, one can find a unique $u\in W$ such that 
\begin{equation} \label{2.22}
\left\{ \begin{array}{ll} L u = 0 & \text{ in } \Omega  \\ 
\Tr u = g & \text{ on } \Gamma,
\end{array}\right.
\end{equation}
where $Lu = 0$ is 
taken in the weak sense \eqref{2.21}. 
The solvability of \eqref{2.22} for any ``smooth enough'' function $g$ proves that the boundary 
is seen by the operator, and so one can hope 
to get analogues of \eqref{2.1}--\eqref{2.2} 
for domains with higher co-dimensional boundaries. 

Starting with the analogue of \eqref{2.1}, we want to show that if $g\in H$ is continuous, then the solution $u\in W$ given by \eqref{2.22} is actually continuous on $\R^n = \overline{\Omega}$. 
The continuity -- actually, the H\"older continuity -- of the weak solutions to $Lu=0$ 
{\em inside the domain} is an immediate consequence of the classical elliptic theory, 
since the weight $w(X)$ is locally bounded from above and below; in a similar spirit, we have  
Harnack estimates inside the domain. The estimates {\em at the boundary} are proven via 
one of the usual methods 
(see for instance \cite[Chapter 8]{GT}): we establish a Cacciopoli inequality {at the boundary} 
and we use Moser's iterations to get Moser estimates. Then 
we prove that a positive solution $u$ to $Lu=0$ with $\Tr u = 1$ on a ball $B$ centered on $\Gamma$ is uniformly bounded from below on $B/2$; for this, we use the Poincar\'e inequality \eqref{2.10} -- which is a way to say $\Gamma$ is ``massive'' enough from the viewpoint of the weighted space $W$ -- 
instead of the classical argument relying on the fatness of the complement. 
In the end, we obtain the following result \cite[Lemma 8.106]{DFMprelim}: 

\begin{theorem}[De Giorgi-Nash-Moser estimates at the boundary] \label{2.23}
Let $B=B(x,r)$ be a ball centered on $\Gamma$ and $u\in W \subset L^2_{loc}(\R^n,dm)$ 
be a weak solution to $Lu=0$ in $\Omega$ such that $\Tr u$ is continuous and bounded on $B$. 
Denote by $\osc_E u$
the difference between the (essential) supremum and the (essential) infimum of $u$ on $E$.
There exists $\alpha >0$ such that for $0<s<r$,
\begin{equation} \label{2.24}
\osc_{B(x,s)} u \leq C \left(\frac{s}{r}\right)^\alpha \osc_{B(x,r)} u 
+ C\osc_{B(x,\sqrt{sr}) \cap \Gamma} \Tr u.
\end{equation}
In particular, $u$ is continuous on $B$. In addition, if $\Tr u = 0$ on $B$, then 
for $x\in \frac12 B$ and $0<s<r/3$, 
\begin{equation} \label{2.25}
\osc_{B(x,s)} u \leq C \left(\frac{s}{r}\right)^\alpha \left( \fint_B |u|^2 \, dm \right)^\frac12 < +\infty.
\end{equation}
The constants $\alpha,C$ depend only on the dimensions $d$ and $n$, as well as the constants $C_0$ and $C_3$.
\end{theorem}

In the previous theorem, we do not need the fact that $u$ is a weak solution to $Lu=0$ everywhere 
in $\Omega$,  but only
in $B$. However, we need to assume {\em a priori} that $u$ is in an appropriate space 
(weaker than $W$, the assumption we made here), but let us 
suppress these technicalities as they are are not essential for construction described below.

In order to construct the elliptic measure, we also need a maximum principle. If $u\in W$ is a weak solution to $Lu=0$ in $\Omega$, then
\begin{equation} \label{2.26}
\sup_\Omega u \leq \sup_\Gamma \Tr u \quad  \text{ and } \quad \inf_\Omega u \geq \inf_\Gamma \Tr u.
\end{equation}
This
may be a surprise at first, since our 
domain $\Omega$ is unbounded and we succeed nevertheless to bound the function $u$ by its trace. 
But we assume here
that $u\in W$, an assumption which therefore cannot be weakened too much, and this 
prevents $u$ from going too
wild at $\infty$.

The solvability of \eqref{2.22} and Theorem \ref{2.23} give us a linear functional $U: \, H \cap C^0_b(\Gamma) \to C^0(\R^n)$, where $C^0_b(\Gamma)$ denotes the space of bounded continuous functions on $\Gamma$ and $C^0(\R^n)$ is the space of continuous functions on $\R^n$, such that $U(g)$ is the unique solution to $Lu=0$ in $W$ satisfying $\Tr u = g$. 
The maximum principle entails that $U$ is a bounded (with the $L^\infty$ norms). 
The Riesz representation theorem for measures then proves the existence of a family of 
Borel regular probability measures $\{\omega^X_L\}_{X\in \Omega}$ on $\dr \Omega$ such that 
\begin{equation} \label{2.27}
U(g)(X) = \int_{\dr \Omega} g(y) \, d\omega^X_L(y) \quad \text{ for } g\in H \cap C^0(\Gamma).
\end{equation}
For a non trivial Borel set $E$, the function $u_E$ defined on $\Omega$ by 
$u_E(X) = \omega_L^X(E)$ does not always lie in $W$, but
we can still prove that $u_E$ a weak solution, in the sense that it lies in $W^{1,2}_{loc}(\Omega)$ 
and verifies \eqref{2.21}.

Having defined elliptic measures $\omega_L$ associated to our degenerate operators $L$, 
the next goal is to show that they have the same suitable  properties mimicking the classical case: non-degeneracy, doubling property, and the comparison principle for harmonic measure. 
To get most of these results, we follow the classical strategy where we build the Green functions and compare them to $\omega_L$.

The Green function is roughly speaking a positive function on $\Omega \times \Omega$ such that, 
for every $Y\in \Omega$, $g(.,Y)$ solves $Lg(.,Y) = \delta_Y$ and $\Tr g(.,Y) = 0$. 
Here $\delta_Y$ represents the delta distribution supported in $\{Y\}$. 
In more precise terms, the Green function is such that for any $f\in C^\infty_0(\Omega)$, 
the function defined by
\begin{equation} \label{2.28}
u(X) = \int_\Omega g(X,Y) f(Y) \, dY \ \text{ belongs to $W_0$ and satisfies $Lu=f$ in the weak sense.}  
\end{equation}
We show that there is a unique function $g$ (the Green function) such that for each $Y \in \Omega$,
$g(.,Y)$ satisfies \eqref{2.28} and has a vanishing trace on $\Gamma$, 
and for which $g(X,.)$ is continuous on $\R^n \setminus \{X\}$ and locally integrable on $\R^n$. 
The existence of the Green function is established as in \cite{GW} (see also \cite{HoK,DK}), 
by taking the weak limit in some appropriate space $S$ of a sequence $g^\rho(.,Y)$ of functions 
in $W_0$ that verify $Lg^\rho(.,Y) = \delta^\rho_Y$, where $\delta^\rho_Y$ is 
an ``appro\-xi\-mation'' of $\delta_Y$ in $(W_0)^*$. 
The main difference with the classical case is that we have two different behaviors when 
$|X-Y| \leq \delta(Y)/2$ and $|X-Y| \geq \delta(Y)/4$, which leads to the proof of two sets 
of estimates, and also the space $S$ in which we take the limit needs to be adapted to this dual nature. 
We also prove that
\begin{equation} \label{2.29}
0 \leq g(X,Y) \leq C|X-Y|^{1-d} \approx C\frac{|X-Y|^2}{m(B(X,|X-Y|))} \quad \text{ if } |X-Y| \geq \delta(Y)/4,
\end{equation}
\begin{equation} \label{2.30}
g(X,Y) \approx \frac{|X-Y|^{2-n}}{w(Y)} \approx \frac{|X-Y|^2}{m(B(X,|X-Y|))} \quad \text{ if } |X-Y| \leq \delta(Y)/2, \, n\geq 3,
\end{equation}
\begin{equation} \label{2.31}
g(X,Y) \approx \frac{1}{w(Y)} \ln\left( \frac{\delta(Y)}{|X-Y|} \right) \quad \text{ if } |X-Y| \leq \delta(Y)/2, \, n=2,
\end{equation}
where the constants above depend only on $d$, $n$, $C_0$, and $C_3$; see \cite{DFMprelim}
for the case when $n \geq 3$ and a slightly weaker estimate when $n=2$, and 
\cite{DFMprelim2} for the remaining case.
As a general rule, the constants that will appear in the rest of the section will have the same dependence
on $d$, $n$, $C_0$, and $C_3$, which will no longer be recalled.

Let us return to the elliptic measure. We want to prove properties similar to the ones found in 
\cite[Section 1.3]{KenigB} in the co-dimension 1 case. 
First, our elliptic measure in higher codimension is non-degenerate, in that for any ball $B$ 
centered on $\Gamma$, one has 
\begin{equation} \label{2.32}
\omega^X(B\cap \Gamma) \geq C^{-1} \ \text{ for } X \in \Omega \cap \frac12 B, \quad  \text{ and } \quad  \omega^Y(\Gamma \setminus B) \geq C^{-1} \quad \text{ for } Y \in \Omega \setminus 2B.
 \end{equation} 
The next property is the comparison between the elliptic measure and the Green function. 
Let $B:= B(x,r)$ is a ball centered on $\Gamma$, and write $A_B \in \Omega$ for 
any corkscrew point associated to $x\in \Gamma$ and $r>0$ given by \eqref{2.7}; then
 \begin{equation} \label{2.33}
\omega^X(\Gamma \cap B) \approx r^{1-d} g(X,A_B)  \ \text{ for } X \in \Omega \setminus  2B,
 \end{equation} 
  \begin{equation} \label{2.34}
\omega^Y(\Gamma \setminus B) \approx r^{1-d} 
g(Y,A_{C_1B}) \quad \text{ for } Y \in \Omega \cap \frac12 B,
 \end{equation} 
where $C_1$ is the constant in \eqref{2.7}. 
The two estimates \eqref{2.33}--\eqref{2.34} can be seen as weak versions of the comparison principle: it is indeed a comparison principle that only deals with Green functions and elliptic measures. 
From \eqref{2.33}--\eqref{2.34}, we deduce the following two key results: the doubling property of 
the elliptic measure, which guarantees that, if $B$ is a ball centered on $\Gamma$, then
\begin{equation} \label{2.35}
\omega^X(2B\cap \Gamma) \leq C \omega^X(B\cap \Gamma) \ \text{ for } X \in \Omega \setminus 4B,
\end{equation} 
and the comparison principle for elliptic measure (also called change-of-pole estimates), which states that
\begin{equation} \label{2.36}
\frac{\omega^X(E)}{\omega^X(\Gamma \cap B)} \approx \omega^{A_B}(E)
\end{equation} 
whenever $B$ is a ball centered on $\Gamma$, $E \subset B\cap \Gamma$ is a Borel set, 
and $X\in \Omega \setminus 2B$. 

More generally, we can deduce 
from \eqref{2.33}--\eqref{2.34} a general version of the comparison principle 
(see \cite{CFMS} in the co-dimension 1 case): if $B$ is a ball centered on $\Gamma$  
and $u,v$ are two positive weak solutions to $Lu=Lv=0$ in $2B$ such that $u,v$ lie in a restriction 
of $W$ to $2B$ and $\Tr u = \Tr v = 0$ on $2B$, then 
\begin{equation} \label{2.37}
\sup_{Z\in B\setminus \Gamma} \frac{u(Z)}{v(Z)} \leq C \inf_{Z\in B\setminus \Gamma} \frac{u(Z)}{v(Z)}.
\end{equation} 
The argument to get \eqref{2.37} is actually a bit more involved than for its analogue in co-dimension 1.
Indeed, in the classical theory, the lack of information on $u,v$ outside $2B$ is compensated 
by comparing the solutions $u,v$ to the elliptic measure associated to the subdomain $\Omega \cap 2B$. 
In the higher co-dimension case, we did not construct an elliptic measure for such domains $\Omega \cap 2B$ whose boundaries have now mixed dimensions, and we compensate this lack of definition  by using the non-degeneracy \eqref{2.35} to transfer estimates on $\dr (2B)$ back to $\Gamma$.

\section{The solvability of the Dirichlet problem in $\Omega_0 = \R^{n}\setminus \R^d$} 
\label{S3}

\subsection{Presentation of the Dirichlet problem with data in $L^p$}
\label{S3.1}
Our further investigation is motivated by the effort to understand the necessary and sufficient geometric conditions for absolute continuity of the elliptic measure with respect to the Hausdorff measure on the boundary. At this point, we start with proving this property on suitable Lipschitz domains. Following the ideas of C. Kenig and his collaborators in the case of co-dimension one, the plan is to prove that 
the harmonic measure is absolutely continuous with respect to the Hausdorff measure if and only if the Dirichlet problem with data in $L^p$ can be solved for some $p\in (1,+\infty)$ and then to establish solvability of the Dirichlet problem by making a change of variables to pass to the ``half-space" (a complement of $\RR^d$ in our case). The problem is that the change of variables affects, of course, the coefficients of the emerging operator $L=-\div A \nabla$  on $\RR^n\setminus \RR^d$ and the goal of the game is to find such a transformation that one can, indeed, prove absolute continuity for the resulting $\omega_L$. The two changes of variables used in co-dimension one, yielding a $t$-independent $A$ and yielding an $A$ satisfying the Carleson measure condition \eqref{CMcoeff} turn out to be not very useful in our case, roughly speaking, because we have many $t$-directions which have to be controlled much more carefully than by a simple analogue of \eqref{CMcoeff}: indeed, we virtually need an isometry in $t$. This seemingly innocent problem yields a completely different construction which is new even in the co-dimension one case and which is the technical core of our argument. Moreover, it turns out that the ``obvious" candidate for the ``Laplacian", $- \div \dist(X,\Gamma)^{d+1-n} \,\nabla$ is not suitable either, and we have to define a new distance function which eventually led us to some mysterious properties (see Section~\ref{S6}). Chronologically, we first established the $A^\infty$ property for a newly found nice operator on lower dimensional Lipschitz graphs and then treated the Dirichlet problem for $1<p<\infty$. We will do the opposite in this survey in order to streamline the exposition and present a slightly different point of view.

Roughly speaking, given a domain  $\Omega \subset \R^n$ and an elliptic operator $L=-\div A  \nabla$ on $\Omega$, the Dirichlet problem is well-posed in $L^p$, $0<p<\infty$, if for every $g\in L^p(\dr \Omega)$ there exists a unique extension $u:=u_g$ that satisfies \eqref{2.1} and suitable estimates on solutions in terms of the data. Let us make this precise. The first line in \eqref{2.1} can be easily interpreted in the weak sense, but the second line of \eqref{2.1}, namely $u = g$ on $\dr \Omega$, is quite unclear at this point if $g$ is neither continuous nor in $H$. 
In addition, we have specified neither the appropriate bounds on solutions in terms of the data, nor the class of functions in which uniqueness will be postulated. 

Analogously to the classical case, the appropriate estimates should be stated as 
$$\|N(u_g)\|_{L^p(\dr \Omega)} \leq C \|g\|_{L^p(\dr \Omega)},$$ with the constant $C>0$  independent of $g$, and the non-tangential maximal function of $u$ is defined for $x\in \dr \Omega$ by 
\begin{equation} \label{3.1}
N(u)(x) = \sup_{X \in \gamma(x)} |u(X)|\qquad \text{ with } \gamma(x) 
= \{X\in \Omega, \, |X-x| \leq C_4
\dist(X,\dr\Omega)\},
\end{equation} 
where we choose $C_4$ larger than $C_1$, so that at least the corkscrew points of \eqref{2.7} 
lie in $\gamma(x)$. The condition $u=g$ on $\dr \Omega$ is then taken as a non-tangential limit; 
more precisely, we ask that
\begin{equation} \label{3.2}
\ds \lim_{X \to x \atop X\in \gamma(x)}
u_g(X) \text{ exists and equals }  g(x) \text{ for a.e. } x\in \dr \Omega.
\end{equation}
Observe that the definition \eqref{3.2} is quite weak but 
has a great synergy with the non-tangential maximum function $N$. Indeed, thanks to the dominated convergence theorem, \eqref{3.2} combined to the assumption that $\|N(u_g)\|_{L^p(\dr \Omega)}$ is finite implies the stronger convergence
\begin{equation} \label{3.3}
 \lim_{r\to 0} \sup_{X \in \gamma(x) \atop \dist(X,\dr \Omega) \leq r}  |u_g(X)-g(x)| = 0 \ \text{ in } L^p(\dr \Omega).
\end{equation}

As mentioned above, using a new change of variables, to be presented in Section \ref{S4}, we reduce the Dirichlet problem on a Lipschitz domain to that in $\RR^n\setminus \RR^d$ for a suitable class of operators.  We expect \eqref{Dp} to be solvable for $L = - \div |t|^{d+1-n} \nabla$ because, as explained in  Section \ref{S2.1}, the solutions to \eqref{Dp} can be obtained by rotating the solutions to $-\Delta u = 0$ in $\R^{d+1}_+$. Indeed, the existence of (radial) solutions to \eqref{Dp} follows from solvability of the Dirichlet problem in $L^p$ for the operator $-\Delta$ and the domain $\R^{d+1}_+$. The uniqueness will be proved in theorems that we state below and will ensure, in particular, that all solutions to this particular problem are radial.
However, much more generally, we will have to address a class of operators which is, in part, dictated by the change of variables (at the very least, it should allow us to treat the Dirichlet problem on a Lipschitz graph for one of our preferred operators), in part by the existing results for co-dimension one and is, in some sense, close to optimal in the sense that at the regularity level, 
 it reaches all the way to existing counterexamples.  The emerging operators and solutions are not necessarily radially independent and no dimension reduction argument is possible.

Returning to the sharpness of our results, some counterexamples from the co-dimension 1 can be adapted to our higher codimension setting. Given an elliptic operator $\cL=-\div \A \nabla$ defined on $\R^{d+1}_+$, we can construct a degenerate elliptic operator $L=-\div A \nabla$ (see a few lines below) such that if $v$ is a solution to $\cL v = 0$ in $\R^{d+1}_+$, then the function $u$ defined on $\Omega_0 = \R^n \setminus \R^d$ by $u(x,t) = v(x,|t|)$ is a solution to $Lu = 0$ on $\Omega$. Moreover, it is easy to check that $u \in W$ if and only if $v\in W_{\R^{d+1}_+} := \{w \in L^1_{loc}(\R^{d+1}_+), \, \nabla w \in L^2(\R^{d+1}_+)\}$; but since the space of traces $H$ is the same for both spaces, the uniqueness in the Lax-Milgram theorem gives that the solution to $Lu=0$ in $\Omega_0$ that lies in $W$ has to be the ones obtained by rotating the solutions in $W_{R^{d+1}_+}$ to $\cL v = 0$. Remember now that the harmonic measure is constructed only from the solutions in $W$, which implies that for $E\subset \R^d$ a Borel set
\begin{equation}
\omega^{(x,t)}_L(E) = \omega^{(x,|t|)}_{\cL}(E).
\end{equation}
So if the harmonic measure $\omega^{(x,|t|)}_{\cL}$ happens to be singular with respect to $\cL^d_{\R^d}$, the Lebesgue measure on $\R^d$, as in \cite{CFK,MM}, 
then the harmonic measure $\omega^{(x,t)}_{L}$ will be of course also singular with respect to $\cL^d_{\R^d}$.

We give now the construction of $A$ (and so $L$) from $\A$. We write $(a_{ij})_{1\leq i,j\leq d+1}$ for the coefficients of the matrix $\A$ and we define
\begin{equation}
A= A(x,t) := |t|^{d+1-n} \begin{pmatrix} A_1 & A_2 \\  A_3 & A_4\end{pmatrix}
\end{equation}
where $A_1$ is the $d\times d$-matrix with coefficients $(a_{ij})_{1\leq i,j\leq d}$, $A_2$ is the $d\times (n-d)$-matrix with coefficients $(\frac{t_j}{|t|}a_{i(d+1)})_{1\leq i \leq d < j\leq n}$, $A_3$ is the $(n-d)\times d$-matrix with coefficients $(\frac{t_i}{|t|}a_{(d+1)j})_{1\leq j \leq d <i \leq n}$, and $A_4 = a_{(d+1)(d+1)} {\rm Id}$, where ${\rm Id}$ is the identity matrix of order $d$. For short, 
\begin{equation}
A:= |t|^{d+1-n} \begin{pmatrix} a_{ij} & \ds \frac{t_j}{|t|} a_{i(d+1)} \\  \ds \frac{t_i}{|t|} a_{(d+1)j} & a_{(d+1)(d+1)} \,{\rm Id} \end{pmatrix}.
\end{equation}
We let the reader check that $L$ is a degenerate elliptic operator satisfying \eqref{2.5}--\eqref{2.6}, and that the claimed properties of $L$ hold.

For the rest of the section, $d<n-1$ is a integer, and we stick to the model case when
$\Gamma = \Gamma_0 := \big\{(x,0)\in \R^n, \, x\in \R^d \big\}$ (which we identify with $\R^d$)
and $\Omega = \Omega_0 := \R^n \setminus \Gamma_0$.
A point of $\R^n$ is often written $X=(x,t)$ or $Y=(y,s)$, where the first coordinate is in $\R^d$ 
and the second one in $\R^{n-d}$, and the coordinates of $t$ are $(t_{d+1}, \dots, t_n)$. 
By Section \ref{S2}, the adequate weight in this case is $w(X) = |t|^{d+1-n}$. 
We address the following problem: given a $p>1$, describe the class of 
degenerate elliptic operators $L:=-\div A  \nabla$, where $A$ satisfies \eqref{2.5}--\eqref{2.6}, 
such that for any $g\in L^p(\R^d)$, one can find a unique $u\in W^{1,2}_{loc}(\R^n\setminus \R^d)$ 
such that
\begin{equation} \label{Dp} \tag{D$_p$}
\left\{\begin{array}{ll}
Lu = 0 & \text{ in } \R^n \setminus \R^d, \\
u = g & \text{ on } \R^d, \\
\|N(u)\|_p \leq C \|g\|_p.
\end{array} \right.
\end{equation}
In \eqref{Dp}, the non-tangential maximal function $N(u)$ is defined for $x\in \R^d$ by \eqref{3.1}, and we  shall work with
$\gamma(x) = \big\{(y,s) \in \Omega, \, |y-x| \leq |s| \big\}$, that is a constant $C_4 = \sqrt 2$ in \eqref{3.1},
but this choice is of little importance.

\subsection{Solvability of the Dirichlet problem in $L^p$ for small $p$} To tackle the well-posedness of the Dirichlet problem, much as in \cite{KKPT,KePiDrift}, and most relevantly for this discussion, \cite{DP},  we want to bound the non-tangential maximal function $N(u)$  with a square function $S(u)$. But in order to  
deal with \eqref{Dp} for $p\neq 2$, we shall replace the square function $S(u)$ with
the $p$-adapted square function $S_p(u)$  introduced in \cite{DPP} and
defined by  
\begin{equation} \label{3.5}
S_p(u)(x) := \left( \iint_{\gamma(x)} |\nabla u|^2 |u|^{p-2} \frac{ds\, dy}{|s|^{n-2}} \right)^\frac1p
\end{equation}
for $x\in \R^d$. The function $S_p(u)$ is defined for all $p\in (1,+\infty)$, and coincides with 
the classical square functional when $p=2$. When $p<2$, multiplying by $|u|^{p-2}$ in the definition of $S_p(u)$ looks problematic because $u$ may vanish on a large set, but observe that for
$u\in W^{1,2}_{loc}(\R^n\setminus \R^d)$,  
$\nabla u \equiv 0$ almost everywhere on $\{u= 0\}$, so
$S_p(u)$ is equal to its variant 
where we only integrate on $\gamma(x) \cap \{\nabla u \neq 0\}$.

The quantity $|\nabla u|^2 |u|^{p-2}$ does not look easy to handle, however, we can link it 
to the notion of $L^p$ dissipativity (see \cite{CD,CM}), and from there, Dindo\v s and Pipher obtained in \cite{DP} (see \cite{FMZ} for the improved version stated here) the following remarkable properties that are the core of our computations, and which hold for any $p\in (1,+\infty)$. First, one has 
a property of $p$-ellipticity, namely
 \begin{equation} \label{3.6}
\iint_{\Omega_0} 
A \nabla u \cdot \nabla [|u|^{p-2}u] \, \chi 
\, dx\, dt \geq C^{-1}_p \iint_{\Omega_0}  |\nabla u|^2 |u|^{p-2} \, \chi 
\, \frac{dt\, dx}{|t|^{n-d-1}}
\end{equation}
whenever $\chi$ is a non-negative bounded function, 
$v:=|u|^{p/2-1}u \in W^{1,2}_{loc}(\Omega_0)$, 
and $\chi \nabla v \in L^2(\Omega_0, \frac{dt\, dx}{|t|^{n-d-1}})$. 
In addition, we can prove the Cacciopoli-type inequality
 \begin{equation} \label{3.7}
  \iint_B  |\nabla u|^2 |u|^{p-2}\, \frac{dt\, dx}{|t|^{n-d-1}} 
  \leq C_p   \iint_B  |u|^{p}\, \frac{dt\, dx}{|t|^{n-d-1}}
\end{equation}
whenever $u$ is a solution to $Lu=0$ and $B$ is a ball satisfying 
$2B \subset \Omega_0$. .

Returning to connections between the square function and the non-tangential maximal function, the goal is to get, for appropriate elliptic operators $L$, the estimates 
\begin{equation} \label{3.8}
\|N(u)\|_p \leq C \|S_p(u)\|_p \quad \text{ and } 
\quad \|S_p(u)\|_p \leq C \|\Tr u\|_p + \epsilon_L \|N(u)\|_p 
\end{equation}
whenever $u \in W$ is a solution to $Lu=0$ such that $\Tr u \in C^\infty_0(\R^d)$, and where $\epsilon_L$ depends on $L$. Once \eqref{3.8} 
is obtained, provided that the quantity $\|N(u)\|_p$ is finite and $\epsilon_L$ is small enough, we get that $\|N(u)\|_p \leq C \|\Tr u\|_p$ for all solutions with smooth trace, and then we get the existence in the Dirichlet problem by a density argument. The uniqueness is rather 
straightforward.

It is, however, essential to know that $\|N(u)\|_p<+\infty$ for a large class of weak solutions $u \in W$ whose trace lies in $C^\infty_0(\R^d)$. That is why we first prove local results in which all involved quantities are at least trivially finite. If $\ell >0$, $E\subset \R^d$, and $e$ is a non-negative $1$-Lipschitz function on $\R^d$, then for $(x,t) \in \R^n$, we introduce the cut-off functions
\begin{equation} \label{3.9}
\chi_{\ell,E,e}(x,t) := \1_{[0,\ell]}(t) \1_{E}(x) \1_{\{e(y) \leq |s|\}}(x,t).
\end{equation}
We write $\chi_{\ell,E}$ for $\chi_{\ell,E,0}$, and also $\chi_{\ell}$ for $\chi_{\ell,\R^d,0}$. 
Observe that if $\ell>0$, $E=B$ is a ball and $e>0$, then $\chi_{\ell,B,e}$ is compactly supported 
in $\R^n\setminus \R^d$. We introduce the local versions of $S_p(u)$ and $N(u)$ as
\begin{equation} \label{3.11}
S_p(u \big|\chi_{\ell,E,e})(x) 
:= \left( \iint_{\gamma(x)} |\nabla u|^2 |u|^{p-2} \chi_{\ell,E,e}\, \frac{ds\, dy}{|s|^{n-2}} \right)^\frac1p
\end{equation}
and
\begin{equation} \label{3.10}
N(u \big|\chi_{\ell,E,e})(x) := \sup_{(y,s) \in \gamma(x) } |u|\chi_{\ell,E,e}.
\end{equation}

We can check that both quantities are finite when 
$u$ is a weak solution to $Lu=0$ (as defined in \eqref{2.21}), 
thanks to \eqref{3.7} and Moser's estimates. 

Furthermore, if $\epsilon >0$ is given, we say that a function $u$ satisfies \eqref{CM} 
if for any ball $B\subset \R^d$ of radius $\ell$, we have 
\begin{equation} \label{CM} \tag{CM$_\epsilon$}
\int_{x\in B} \int_{|t| \leq \ell} |u|^2 \frac{dt \, dx}{|t|^{n-d}} \leq \epsilon \ell^d.
\end{equation}
If $U=(u_{ij})$ is a matrix, then $U$ satisfies \eqref{CM} if all its coefficients $u_{ij}$ satisfy \eqref{CM}. 
In the co-dimension 1 case, $t$ is positive  
and $|t|^{n-d} = t$.

\begin{theorem}\cite{FMZ} \label{3.12} 
Let $L:= -\div A\nabla$ be an elliptic operator on $\Omega_0 = \R^n \sm \R^d$
satisfying \eqref{2.5}--\eqref{2.6}. 
Assume that the reduced matrix $\A:= |t|^{n-d-1} A$ can be written as the matrix by blocks $[d+(n-d)] \times [d+(n-d)]$
\begin{equation} \label{3.13}
\A := \begin{pmatrix} \A_1 & \A_2 \\ \cB_3+ \C_3 & b\,{\rm {\rm Id}} + \C_4 \end{pmatrix}
\end{equation}
where $C^{-1} \leq b \leq C$ is a scalar function, and $|t||\nabla \cB_3|$, $\C_3$, $\C_4$, and $|t||\nabla b|$ all satisfy \eqref{CM} for some $\epsilon>0$ (that does not need to be small). 
Then for all choices of $p,q\in (1,+\infty)$,  
$\ell>0$, a ball $B\subset \R^d$ with center $x_B$ and radius at least $\ell$,  
a $1$-Lipschitz function $e>0$, and a weak solution $u$ to $Lu=0$, one has
\begin{equation} \label{3.14}
\|S_q(u|\chi_{\ell,B,e})\|_p \leq C \|N(u \big|\chi_{2\ell,2B,e/2})\|_p
\end{equation}
and
\begin{equation} \label{3.15}
\|N(u|\chi_{\ell,B,e})\|_p \leq C \|S_q(u \big|\chi_{2\ell,2B,e/2})\|_p + C\ell^d |u(x_B,\ell)|,
\end{equation}
where $C$ depends only 
on $d$, $n$, $p$, $q$, $\epsilon$, and the ellipticity 
constant $C_3$.
\end{theorem}

The above theorem is a variation of \cite[Theorem 1.7]{FMZ}. 
An important feature of this result is that there is no assumption on the first $d$ lines 
of $A$ (other than the ones coming from ellipticity). Concerning the lower right block, it is 
important for the proof that the main piece $b \,{\rm Id}$ is a (slowly moving) scalar function
times the identity. That is, all the directions in $\R^{n-d}$ need to be treated roughly 
the same way, with errors measured in terms of Carleson measures. In co-dimension $1$,
this issue does not arise, because there is only one coordinate.

As far as we know, \cite{DFMAinfty,FMZ} is the first time where one proves
such local estimates directly (in anterior papers, they would be derived from global ones) and it is a valuable asset particularly for complex coefficient operators when the existing arguments of reduction from global to local bounds fail.

In the co-dimension 1 case, a key argument of the proof of \eqref{3.15} is the fact that we can deduce estimates of the type $\|\Tr u\|_p \leq \|S_2 u\|_p$ for any Lipschitz subdomain (namely sawtooth domains) of $\R^{d+1}_+$ from the bound $\|\Tr u\|_p \leq \|S_2 u\|_p$ when the domain is $\R^{d+1}_+$ by using a simple bi-Lipschitz change of variables that maps the Lipschitz sub-domain back to $\R^d$. 
In our higher codimension case, the analogues of these saw-tooth domains are the sets $\{(x,t)\in \R^n \, |t| > \varphi(x)\}$, where $\varphi:\, \R^d \to [0,+\infty)$ is Lipschitz. These latter sets have mixed co-dimension, depending on whether $\varphi(x) = 0$ or $\varphi >0$, and cannot be mapped back into $\R^n \setminus \R^d$ by a bi-Lipschitz change of variables. Our freedom to choose $e$ as a 1-Lipschitz function (that restrains our estimates inside a saw-tooth domain) solves this difficulty. Actually, the properties of the functions $\chi^1=\chi_{\ell,B,e}$ and $\chi^2=\chi_{2\ell,2B,e/2}$ that really matter are the fact that they are both characteristic functions 
of Lipchitz domains with bounded Lipschitz characters 
and that we can find a smooth function $\psi$ such that $\chi^1 \leq \psi \leq \chi^2$ that satisfies $|\nabla \psi(x,t)| \leq 100/|t|$. 

\medskip

We aim to let $e$ tend to $0$, $B$ tend to $\R^d$, and $\ell$ to $+\infty$ 
in the estimates \eqref{3.14}--\eqref{3.15}. But nothing guarantees the finiteness of the resulting quantities. This  
is why we prove the following result (see \cite[Theorem 1.10]{FMZ}):

\begin{theorem} \label{3.16}
Let $L:= -\div A\nabla$ be an elliptic operator satisfying \eqref{2.5}--\eqref{2.6}. Let $u$ be 
a special solution to $Lu=0$, i.e., such that $u\in W$ 
and $\Tr u \in C^\infty_0(\R^d)$. Then, for any $p\in (1,+\infty)$, one has
\begin{equation} \label{3.17}
\|u\|_{W^p}:= \left(\int_{\R^n\setminus \R^d} |\nabla u|^2 |u|^{p-2} \frac{dt\, dx}{|t|^{n-d-1}}\right)^\frac1p < +\infty.
\end{equation}
\end{theorem}

Observe that $\|.\|_{W^2}$ is simply $\|.\|_{W}$. For any $\ell >0$ and any special  
solution u, notice also that $\|S_p(u|\chi_\ell)\|_p \leq \ell \|u\|_{W^p} < +\infty$. The latter combined to \eqref{3.14}--\eqref{3.15} implies that for such solutions, the quantities $\|S_q(u|\chi_\ell)\|_p$ and $\|N(u|\chi_\ell)\|_p$ are all finite. We then prove that for special  
solutions, and under the assumptions of Theorem \ref{3.12}, one has
\begin{equation} \label{3.18}
\|N(u|\chi_\ell)\|_p \leq C \|S_p(u|\chi_{2\ell})\|_p + C\|\Tr u\|_p,
\end{equation}
and
\begin{equation} \label{3.19}
\|S_p(u|\chi_{2\ell})\|_p \leq C\epsilon \|N(u|\chi_{4\ell})\|_p + C\|\Tr u\|_p + \mathcal E_\ell(u),
\end{equation}
where $\epsilon$ is the Carleson bound on the coefficients of $A$ assumed in Theorem \ref{3.12}, 
and where $\mathcal E_\ell(u)$ is an extra term (which can take any real value) due to the cut-off at the level $\ell$. 
We prove that 
$\|N(u|\chi_\ell)\|_p \leq C\epsilon \|N(u|\chi_{4\ell})\|_p + C\|\Tr u\|_p + \mathcal E_\ell(u)$,
and then that 
$\|N(u|\chi_{\ell_m})\|_p \leq C\|\Tr u\|_p$ for a sequence of $\ell_m$ that goes to $+\infty$,
under the condition that $\epsilon$ is small enough [we use the fact that 
$\|N(u|\chi_{4\ell}) - \|N(u|\chi_\ell)\|_p$ is small and $\mathcal E_\ell(u)$ is non-positive for many $l$ when 
$\|u\|_{W^p}$ is a special 
solution]. Finally, we  prove the existence of a solution $u$ to \eqref{Dp} when $g\in C^\infty_0(\R^d)$ and $L$ satisfies the assumptions of Theorem \ref{3.12} with $\epsilon \leq \epsilon_0(p)$ small enough. From there, it is not very difficult to obtain the following:

\begin{theorem} \label{3.20}
Let $L:= -\div A\nabla$ be an elliptic operator satisfying the same assumptions as in 
Theorem \ref{3.12}. For all $p\in (1,+\infty)$, there exists $\epsilon_0(p)$ such that if $\epsilon \leq \epsilon_0(p)$ -- where $\epsilon$ is the Carleson measure bound on the coefficients assumed in Theorem \ref{3.12} -- then the Dirichlet problem \eqref{Dp} is solvable.
\end{theorem}

Remarkably, the proof of Theorem \eqref{3.20}  can be adapted to the case 
where $A$ has
complex coefficients. The estimates \eqref{3.6}--\eqref{3.7} do not hold for the full range of $p$ in $(1,+\infty)$ anymore, but only on a restricted range $(p_0,p_0')$ around $2$, and the solvability of the Dirichlet problem \eqref{Dp} will be obtained for the same range of $p$ as long as we change the definition of non-tangential maximal function $N$ to an averaged version, often denoted by $\wt N$.

\subsection{Solvability of the Dirichlet problem in $L^p$ for large $p$ and in BMO. Absolute continuity of the elliptic measure with respect to the Hausdorff measure on the boundary}
\label{S3.3}

While the smallness of $\eps>0$ in Theorem~\ref{3.20} is necessary to establish well-posedness of \eqref{Dp} for a given $p$ (particularly when $p$ is small), we show that the mere finiteness of the involved Carleson norms of coefficients is necessary for the well-posedness of \eqref{Dp} for {\it some} $p<\infty$. This is very important due to the aforementioned fact that the solvability of \eqref{Dp} for some $p<+\infty$ is actually equivalent to the mutual absolute continuity of the elliptic measure (associated to $L$, as defined in Section \ref{S2}) and the $d$-dimensional Hausdorff measure $\sigma$, 
and more precisely to the fact that $\omega_L \in A_\infty(\sigma)$, with the ad hoc definition given 
near \eqref{3.23}. 
In our setting the proof of this fact is given in \cite[Theorem 4.1]{MZ}.

Concerning $A_\infty$, we were able to obtain (see \cite[Theorem 1.32]{DFMAinfty}) the following result.

\begin{theorem} \label{3.22}
Let $L:= -\div A\nabla$ be an elliptic operator on $\Omega_0 = \R^n \sm \R^d$ 
that satisfies the assumptions of
Theorem \ref{3.12}. Then $\omega_L \in A^\infty(\sigma)$ (as defined near \eqref{3.23}).
\end{theorem}

The proof of Theorem \ref{3.22} follows an idea developed in \cite{KKPT}, \cite{KKiPT}, \cite{DPP2} and uses the estimate \eqref{3.14} given above. In short, the idea is that assuming the failure of the $A^\infty$ property, one can construct a boundary set whose characteristic function gives rise to rapidly oscillating solutions (the elliptic measure), and ultimately, to a large square function since the gradient of the solution is systematically large in Whitney balls.

We mention, parenthetically, that in this and a much more general geometric setting (on all AR domains) the property $\omega_L \in A^\infty(\sigma)$ is also equivalent to the BMO solvability of the Dirichlet problem. To be a bit more specific, we need some definitions. 

Let $\Gamma\subset \RR^n$ be a $d$-Ahlfors regular set, where $d<n-1$ is an integer. For any $x\in\Gamma$ and $r>0$, we use $\Delta = \Delta(x,r)$ to denote the surface ball $B(x,r)\cap \Gamma$, and use $T(\Delta) := B(x,r)\cap\Omega$ to denote the ``tent'' above $\Delta$.
A function $f$ defined on $\Gamma$ is a BMO function if
	\begin{equation}\label{eq:defBMO}
		\|f\|_{BMO} := \sup_{\Delta \subset \po} \left(\fint_{\Delta} |f-f_\Delta|^2 d\sigma\right)^{\frac{1}{2}} <\infty,
	\end{equation}
where $f_\Delta$ denotes the average $\fint_{\Delta} f d\sigma$.

We say that the Dirichlet problem \eqref{2.1} is solvable in BMO if  for any boundary function $f\in C_0^0(\po)$, the solution $u$ to \eqref{2.1} given by \eqref{2.27} (suitably extended, see \cite{DFMprelim}) is such that $|\nabla u|^2 \delta(X)^{d-n+2}\,dX$ is a Carleson measure with norm bounded by a constant multiple of $\|f\|_{BMO}^2$, that is, 
	\begin{equation}\label{MT:Carlest}
		\sup_{\Delta\subset\po} \frac{1}{\sigma(\Delta)} \iint_{T(\Delta)} |\nabla u|^2 \delta(X)^{d-n+2}\,dX \leq C \|f\|_{BMO}^2.
	\end{equation}

\begin{theorem}\label{thm:main}\cite{MZ}
	Let $\Gamma$ be a $d$-Ahlfors regular set in $\RR^n$ with $d<n-1$ and $\Omega = \mathbb{R}^n \setminus \Gamma$. Consider the operator $L=-\div  A \nabla$ with a real, symmetric $n\times n$ matrix $A$ satisfying \eqref{2.5} and \eqref{2.6}.
	Then the harmonic measure $\omega \in A_{\infty}(\sigma)$ if and only if the Dirichlet problem is solvable in BMO.
\end{theorem}
In co-dimension 1 this has been proved in \cite{DKP} for Lipschitz domains and in \cite{Zh} for uniform domains with Ahlfors regular boundaries. In the present setting the result is in \cite{MZ}. The reader can see that, in particular,  it entails BMO solvability of the Dirichlet problem for the operators from Theorem~\ref{3.22} and the upcoming Theorem~\ref{4.9}.

\section{The case when
the boundary is a small Lipschitz graph} 
\label{S4}

In this section, we return  
to the discussion of  
Subsection \ref{S2.1}, where we are given the domain $\Omega = \R^n \sm \Gamma$
bounded by an Ahlfors regular set $\Gamma$ of dimension $d < n-1$, and we want to 
see whether we can find degenerate elliptic operators $L$, as in Section \ref{S2},
such that $\omega_L \in A_\infty(\sigma)$, where $\sigma = \H^d_{\vert \Gamma}$
(or any measure on $\Gamma$ satisfying \eqref{defADR}). 
We expect the answer to depend on regularity properties of both $\Gamma$ and the
matrix $A$, but here we shall try to start with Lipschitz graph domains and choose the best operator $L$ in terms of $\Gamma$ so that $\omega_L \in A_\infty(\sigma)$. This was our original incentive and a part of the reason for developments in Section~\ref{S2}.

In the co-dimension 1 case, the simplest operator and the starting point of all investigations is the Laplacian. In the domains with lower dimensional boundaries, there is no clear choice of best operator $L$, but let us decide to
pick $L = - \div[A(X) \nabla]$, with a matrix $A$ which is the product of a scalar positive function
by the identity matrix. This makes sense because we prefer $L$ to be reasonably isotropic, and
taking a complicated matrix that potentially anihilates the effects of a bad changes of variable 
would not seem reasonable. Given the size condition \eqref{2.5}--\eqref{2.6}, the simplest choice 
would seem to be
\begin{equation} \label{Ldist}
L = - \div[\delta(X)^{d+1-n} \nabla] = - \div[\dist(X,\Gamma)^{d+1-n} \nabla],
\end{equation}
but (after an unsuccessful attempt to prove the desired result when $d \geq 2$) we realized that
$\delta(X)$ is not always smooth enough, roughly speaking displaying the same bad features as P. Jones $\beta_\infty$ coefficients in the dimensions larger than 3, and we decided to replace $\delta(X)$ with any of the 
softer functions $D_\alpha(X)$ defined as follow.

Let $\sigma$ be any measure on $\Gamma$ satisfying \eqref{defADR} (but the restriction to $\H^d$
is a good choice), and let $\alpha > 0$ be given. Set $D_\alpha(X) = 0$ for $X\in \Gamma$ and
\begin{equation} \label{4.8}
D_\alpha(X) = \left( \int_\Gamma |X-y|^{-d-\alpha} d\sigma(y) \right)^{-1/\alpha}
\end{equation}
otherwise. It is not hard to deduce from \eqref{defADR} that $D_\alpha(X)$ is equivalent to $\delta(X)$
(and this is the reason for the strange power). The distance $D_\alpha$  has an intrinsic
flavor and  does not see corners of $\Gamma$ as well as the Euclidean $\delta$ does. We set
\begin{equation} \label{LDa}
L_\alpha = - \div[D_\alpha(X)^{d+1-n} \nabla],
\end{equation}
and use this as our (family of) preferred operators.  Both $L_\alpha$ and $D_\alpha$ have some nice and even sometimes mysterious properties as we will describe in Section~\ref{S6}. We will present here some arising questions,  partial answers, and new hopes. For now let us start with the much awaited positive result.

\begin{theorem} \label{4.9}
Let $\varphi:\, \R^d \to \R^{n-d}$ is a Lipschitz function, and let 
$\Gamma = \big\{ (x,\varphi(x)), \, x\in \R^d \big\} \subset \R^n$ denote its graph. 
Set $\sigma = \H^d_{\vert \Gamma}$ and let $\alpha > 0$ be given. 
Define $L_\alpha$ on $\Omega = \R^n \sm \Gamma$ by \eqref{LDa}.
If the Lipschitz constant $\|\varphi\|_{Lip}$ is small enough, depending only on $n$, $d$, and $\alpha$,
then $\omega_{L_\alpha} \in A^\infty(\sigma)$ (as defined near \eqref{3.23}),
with $A^\infty$ constants that depend only on $n$ and $\alpha$. 
\end{theorem}

This is \cite[Theorem 1.18]{DFMAinfty}). It stays true when $\sigma$ is any Ahlfors regular measure
whose support is $\Gamma$ (i.e., such that \eqref{defADR} holds), and then the $A^\infty$ constants
depend also on the constant in \eqref{defADR}. 
It is also true when $d=1$ and $L$ is given by
\eqref{Ldist}, but the proof fails in higher dimensions, due to the lack of regularity of $\delta$. 

This result can be seen as the slightly weaker analogue in higher co-dimensions of Dahlberg 
result of $A_\infty$ absolute continuity of harmonic measure on Lipschitz graphs \cite{Da77}.

The argument described below does not extend to general Lipschitz graphs, but the authors expect
the result to hold for those (with a different proof), and even for uniformly rectifiable sets. However, as this text is written we are still verifying some pertinent details.

Let us emphasize that when $\Gamma$ is a Lipschitz graph, or even the image of
$\Omega_0 = \R^n \sm \R^d$ by any bi-Lipschitz mapping $\rho$, the existence of a degenerate 
elliptic operator $L$ such that $\omega_L\in A^\infty(\sigma)$ is immediate (but not so interesting). Indeed we can take for $L$ the conjugate by $\rho$ of the model operator 
$L_0:= -\div |t|^{d+1-n} \nabla$ on $\Omega_0$, and deduce the result for $L$ from
the result for $L_0$, as suggested in Subsection \ref{S3}. This is why we decided to stick to
simple matrices $A$. In fact, our proof applies to a class of matrices satisfying a certain form of controlled Carleson oscillation condition, but let us not state it here.

We could always hope to prove a converse, saying that if $\omega_{L_\alpha} \in A^\infty(\sigma)$
then $\Gamma$ is uniformly rectifiable, but at this point we do not know how to do that, and, in fact, there are specific values of $\alpha$ for which this fails; see Section \ref{S6}.

Let us now describe the proof of Theorem \ref{4.9}. We start with a few observations.
As mentioned in the introduction, the precise notion of connectedness that is needed for the $A_\infty$ properties has been
causing some trouble in co-dimension $1$. 
Remarkably, it is not an issue here, because \eqref{2.7} and \eqref{cs1}-\eqref{cs2} are automatically verified, since $\Gamma$ is too small to bar the passage. 

The general strategy for the proof is to use a bi-Lipschitz change of variables,
as suggested in Subsection \ref{S3}, to reduce to a degenerate elliptic operator on 
$\Omega_0 = \R^n \sm \R^d$ that we can control, and the matters have been designed so that  Theorem \ref{3.22} gives a sufficiently large class of operators (those described in Theorem \ref{3.12}), so that
the conjugate of $L_\alpha$ by the best change of variables that we could find lies in that class. This change of variables is new even in the classical setting and is an important part of our argument. Let us discuss the details.

Let $\rho: \R^n \to \R^n$ be a bi-Lipschitz change of variables that maps $\Gamma_0 = \R^d$
to $\Gamma$, the graph of $\varphi$. We write $J$ for the Jacobian matrix of $\rho$. 
If $L= - \div A \nabla$ and $u$ is a weak solution to $Lu=0$, one can check that 
$v:= u\circ \rho$ is a weak solution to $L_\rho v:= \div A_\rho  \nabla v = 0$ where
\begin{equation} \label{4.1}
A_\rho(X) = |\det(J(X))| J(X)^{-T} A(\rho(X)) J(X)^{-1}
\end{equation}
and $J^{-T}$ denotes the transpose of the inverse of $J$, 
and as before we shall denote by $\A_\rho$ the reduced matrix $|t|^{n-d-1}A_\rho$. 

Let us review the change of variables used in the past, and at the same time explain some of the
constraints of our setting. The simplest bi-Lipschitz mapping $\rho$ is given by
\begin{equation} \label{4.2}
\rho_1(x,t) = (x,t+\varphi(x)).
\end{equation}
It is easy to see that $\rho_1$ maps $\R^{n} \setminus \R^d$ to $\Omega= \R^n \setminus \Gamma$, 
and that the associated Jacobian matrix $J_1$ is independent of $t$. 
If $L=-\Delta$, or more generally if $A$ is independent of $t$, then $A_{\rho_1}$ also is $t$-independent;
this partially explains the long history of study on elliptic operators with $t$-independent coefficients in the 
co-dimension 1 case. However, in higher co-dimension, the ellipticity conditions \eqref{2.5}--\eqref{2.6} force $A(X)$ to depend on the distance to $\Gamma$, and so the corresponding $A_{\rho_1}(X)$ cannot be $t$-independent. We could also hope to apply Theorem \ref{3.22}, but in this respect the next attempt 
is better.

The second classical change of variables, used for the Dirichlet problem in $L^p$ in \cite{KePiDrift}
(but attributed in different contexts to Necas, Stein, and Kenig) is
\begin{equation} \label{4.3}
\rho_2(x,t) = (x,ct+\eta_{|t|}\ast 
\varphi(x)), \qquad \text{ where $c>0$ and $\eta_r$ is a mollifier.} 
\end{equation}
If the constant $c$ is chosen carefully (depending on the mollifier $\eta_{r}$ and 
the Lipschitz norm $||\varphi||_{lip}$) 
$\rho_2$ is a bi-Lipschitz map from $\Omega_0$ to $\Omega$. 
The Jacobian matrix $J_2$ is such that $|t||\nabla J_2|$ satisfies \eqref{CM} for some $\epsilon>0$ depending on $\eta_{r}$ 
and the Lipschitz constant of $\varphi$. So, as long as the gradient of $A$ satisfies some 
Carleson estimates in the spirit of \eqref{CM} (which is the case for $A=|t|^{d+1-n}{\rm Id}$), 
then $|t| |\nabla \A_{\rho_2}|$ shall verify (CM$_{\epsilon'}$). 
This change of variables suffices in co-dimension 1, where saying that $|t| |\nabla \A_{\rho_2}|$ satisfies (CM$_{\epsilon'}$), e.g.,  fits into the profile \eqref{3.13}. However, it is not suitable in our higher co-dimension setting, because the bottom right corner of $\A_{\rho_2}$ cannot always 
be written as $b\,{\rm Id}+\C_4$ as requested in \eqref{3.13} (at least as long as $A$ is written as a scalar function times the identity).

The conclusion is that the two previously used changes of variable do not help in our setting. 
We observed that the choice of $\rho_2$ didn't work because $\rho_2$ fails to conserve 
the scalar form of the bottom right corner of $A$. Therefore, we try to construct our change 
of variable $\rho$ so that it is nearly conformal in the $t$-coordinates. We take 
\begin{equation} \label{4.4}
\rho(x,0) = (x,\varphi(x)) \quad \text{ for } x\in \R^d
\end{equation}
and
\begin{equation} \label{4.5}
\rho(x,t) = (x,\eta_{|t|}*\varphi(x)) + h(x,|t|) \, R_{x,|t|}(0,t) \quad \text{ for } (x,t)\in \Omega_0,
\end{equation}
where $R_{x,r}$ is a linear isometry of $\R^n$ and $h(x,r) > 0$ is a dilation factor,
that we rapidly discuss now. We construct $R_{x,r}$ (with a convolution formula and projections)
so that it maps $\R^d$ to the $d$-plane $P(x,r)$ tangent to 
$\Gamma_r := \{(x,\eta_r*\varphi(x)), \, x\in \R^d\}$ at the point 
$\Phi_r(x):= (x,\eta_r*\varphi(x))$, and hence also $R_{x,r}$ maps $\R^{n-d} = (\R^d)^\bot$ 
to the orthogonal plane to $P(x,r)$ at $\Phi_r(x)$. The additional dilation factor $h(x,r)$ is perhaps
not even needed in the present case, but it gives a little bit of extra flexibility and helps checking 
that when we do the composition by $\rho$, we get a matrix $A_\rho$ that satisfies the 
assumptions of Theorem \ref{3.12}. We need $\|\varphi\|_{lip}$ to be small enough, because we can
only prove that our formula gives a bi-Lipschitz mapping under this condition.

Let us just say a few words about the control on $A_\rho$.
Analogously to the change of variables $\rho_2$, the conjugate $A_\rho$ of the matrix $A$ by $\rho$ 
will be of the form \eqref{3.13} if, roughly speaking, $A$ satisfies a condition similar to \eqref{3.13} 
in the first place. Here we took $A= g \, {\rm Id}$, where $g$ is a scalar function equivalent to the weight $w$. 
In this case, we can show that $A_\rho$ satisfies 
\eqref{3.13} if we can find a decomposition 
\begin{equation} \label{4.6} 
\frac{g \circ \rho(x,t)}{|t|^{d+1-n}} = 1 + f, \quad \text{ where $f$ satisfies (CM$_\epsilon$) for some $\epsilon>0$.}
\end{equation}
To this end some regularity for the weight $g$ is helpful.
The choice of $g(X) = w(X) = \delta(X)^{d+1-n} =  \dist(X,\Gamma)^{d+1-n}$ works when 
$d=1$, and this has to do with the fact that we can 
link the distance $\dist(X,\Gamma)$ to the Peter-Jones $\beta$-number 
\begin{equation} \label{4.7}
\beta_\infty(x,r) = \inf_{a \text{ affine}} \sup_{y\in B(x,r)} |\varphi(y) - a(y)|,
\end{equation}
and that 
$|\beta_\infty(x,t)|^2 dx\, dt/t$ is a Carleson measure when $d=1$.
As a result, we can show that \eqref{4.6} is true when $g= w$ and $d=1$.
When $d>1$, the weight $w$ {\em a priori} lacks smoothness, \eqref{4.6} may well fail, 
and we choose the smoother weight $g = D_\alpha^{d+1-n}$, where $D_\alpha$
is as in \eqref{4.8}. The difference with $\delta(X)$ is that the less local nature of
$D_\alpha$ allows us to estimate its variations in terms of Tolsa's $\alpha$-numbers, 
whose definition is omitted here (but see \cite{Tolsa09}), and then use the fact that for uniformly
rectifiable sets -- hence also Lipschitz graphs -- $|\alpha(x,t)|^2 dx\, dt/t$ is a Carleson measure,
to prove that $g := D_\alpha^{d+1-n}$ satisfies \eqref{4.6}.

\section{The miraculous case: explicit Green function} \label{S6} 
One of the first points that we teach in a PDE class is that there are very few known explicit solutions. In the classical co-dimension one case one can write a Green function for a half space, a ball, and maybe for very few polygons. It came as a big surprise to the first and the last authors of this paper and M. Engelstein \cite{DEM} that in the higher co-dimensional case there is a choice of $\alpha>0$ such that an explicit solution (and in most cases even a Green function with the pole at $\infty$) for $L_\alpha$ from \eqref{LDa} can be written explicitly on a domain with {\it any} AR boundary.  

If the numbers $n$, $d < n$, and $\alpha > 0$ are such that
\begin{equation}\label{5.1}
n = d+2+\alpha
\end{equation}
then it turns out that the function $D_\alpha$ defined in \eqref{4.8} is also a solution of $L_\alpha D_\alpha = 0$ in $\Omega=\RR^n\setminus \Gamma$,
where $L_\alpha$ is the degenerate elliptic operator that we like 
to associate to $E$, \eqref{LDa}. The reader can check this by an explicit computation and the equality even holds strongly in $\Omega$ since $D_\alpha$ is smooth.

Notice that there are absolutely no restrictions on the underlying AR set $\Gamma$, it could even have a fractional dimension.

Morally speaking, this gives the Green function with a pole at infinity: indeed, $D_\alpha$ is a solution vanishing at the boundary and not growing too fast as $X$ goes to infinity. In reality, the questions of uniqueness, boundary limits, etc.  become more involved in such general domains, but at the very least 
we get the (surprising!) fact that the harmonic measure
$\omega^X_{L_\alpha}$ is not only (mutually) absolutely continuous to $\sigma$, but even comparable, as follows.

\begin{theorem}\cite{DEM}
Let $\Gamma \subset \R^n$ be an Ahlfors regular set of dimension $d<n-2$ and $\alpha=n-d-2$. Let $X \in \Omega = \R^n \sm \Gamma$ be given, and denote by
$\omega^X$ the associated harmonic measure with pole at $X$. Set $R = \dist(X,E)$.
Then there is a constant $C$, that depends only on $n$, $d$, and the Ahlfors regularity constant 
for $\Gamma$, such that
\begin{equation} \label{5.4}
C^{-1} \sigma(A) \leq R^d \omega_{L_\alpha}^X(A) \leq  C \sigma(A)
\ \text{ for every measurable set } A \subset \Gamma \cap B(X,100R).
\end{equation}
Here $L_\alpha$ is the operator given by \eqref{LDa}.
\end{theorem}

It follows that in the special case when $n = d+2+\alpha$, there is absolutely no converse to the theorem that
says that $\omega^X \in A^\infty(\sigma)$ when $E$ is uniformly rectifiable. Even sets of fractional dimensions work! This is, of course, specific to the particular choice of $\alpha$ and a particular choice of the operator, and it is entirely possible that for all other operators of our interest the converse statement is indeed valid.

\end{document}